# INTERMEDIATE OPTIMAL GEVREY EXPONENTS OCCUR

MICHAEL CHRIST

ABSTRACT. Hypoellipticity in Gevrey classes $G^s$ is characterized for a simple family of sums of squares of vector fields satisfying the bracket hypothesis, with analytic coefficients. It is shown that hypoellipticity holds if and only if $s$ is greater than or equal to an optimal exponent that may take on any rational value.

## 1. INTRODUCTION

Let $L$ be a linear partial differential operator, with $C^\infty$ coefficients. Denote by $G^s$ the Gevrey class of order $s$, for each $1 \le s < \infty$. The operator $L$ is said to be $G^s$ hypoelliptic at a point $x$ if, whenever $u$ is a distribution defined near $x$ for which $Lu \in G^s$ in some neighborhood of $x$, necessarily $u \in G^s$ in some neighborhood, possibly smaller, of $x$. It is $G^s$ hypoelliptic in an open set $U$ if it is so at every point of $U$.

**Definition.** The optimal exponent for Gevrey regularity for $L$ at $x$ is defined to be the infimum of all $s \ge 1$ such that $L$ is $G^s$ hypoelliptic at $x$. The optimal Gevrey exponent for $L$ in an open set $U$ is defined to be the supremum, over all $x \in U$, of the optimal exponent for $L$ at $x$.

Suppose that $L$ is locally solvable in $L^2$ near $x$, that is, that there exists a neighborhood $V$ of $x$ such that $\|\varphi\|_{L^2} \le C\|L^*\varphi\|_{L^2}$ for all $\varphi \in C_0^\infty(V)$, where $L^*$ denotes the transpose of $L$. Under this mild hypothesis, if $L$ is $G^s$ hypoelliptic in $V$ for some exponent $s$, then it is also $G^t$ hypoelliptic in $V$ for every $t \ge s$; see Theorem 3.1 and Remark 3.2 of Métivier [9]. The word "optimal" is thus appropriate under this hypothesis.

Suppose that $L$ is a sum of squares of real vector fields having real analytic coefficients. Fix a point $x$, and denote by $m \ge 1$ the least integer for which the vector fields satisfy the usual bracket hypothesis to order $m$ at $x$. Then in some neighborhood of $x$, $L$ is $G^s$ hypoelliptic for every $s \ge m$ [5]. But $s = m$ need not be optimal; for instance, if $L$ has symplectic characteristic variety and loses one derivative then $L$ is

*Date*: March 8, 1996.

Research supported by NSF grant DMS-9306833 and at MSRI by NSF grant DMS-9022140.

analytic hypoelliptic [11],[13], so $G^s$ hypoelliptic for all $s \geq 1$. Various examples and special classes that are not analytic hypoelliptic have been studied, and an analysis of the proofs that these fail to be analytic hypoelliptic near $x$ has generally shown that they also fail to be $G^s$ hypoelliptic whenever $s < m$. To this author's knowledge, it has only recently been pointed out [3] that there exist sums of squares operators for which the optimal Gevrey exponent satisfies $1 < s_0 < m$, and moreover, $s_0$ need not be an integer but can take on (certain) rational values.

In this note, examples will be given for which $s_0$ takes on intermediate values in $(1, m)$.[1] The analysis of these examples is simpler than that of the operators of [2],[3]. In addition, some general results concerning the connection between Gevrey regularity, variants of the FBI transform, and holomorphic extendibility will be established, generalizing well known results for the real analytic case $s = 1$.

The operators that we study have previously been considered by Grušin [6], who proved certain of them to be analytic hypoelliptic, and by Oleĭnik and Radkevič [10], who proved all others not to be analytic hypoelliptic. More recently they have been rediscovered by Hanges and Himonas, who observed that their principal symbols have symplectic characteristic varieties, even in the analytic nonhypoelliptic case.

Let $1 \leq p \leq q \in \mathbb{N}$, let $(x, t)$ be coordinates in $\mathbb{R} \times \mathbb{R}^2$, and define

$$L = \partial_x^2 + x^{2(p-1)}\partial_{t_1}^2 + x^{2(q-1)}\partial_{t_2}^2. \tag{1.1}$$

These are analytic hypoelliptic if [6] and only if [10] $p = q$. The bracket condition is satisfied to order $m = q$ at $0$.

**Theorem 1.1.** *$L$ is $G^s$ hypoelliptic in some neighborhood of $0$ if and only if $s \geq q/p$.*

Further generalizations are possible. For instance,

$$L = \partial_x^2 + a_1(x, t)x^{2(p-1)}\partial_{t_1}^2 + a_2(x, t)x^{2(q-1)}\partial_{t_2}^2$$

satisfies the same conclusions provided that $a_j$ are real-valued, are strictly positive and are analytic.[2] The variable $(x, t_1, t_2)$ can be permitted to take values in $\mathbb{R}^{m_1} \times \mathbb{R}^{m_2} \times \mathbb{R}^{m_3}$ for any $m_j \geq 1$, with $\partial_x^2, \partial_{t_j}^2$ replaced by constant-coefficient elliptic operators having nonpositive principal symbols and $x^2$ replaced by $|x|^2$.

---

[1] For examples (1.1), $s_0^{-1}$ is always a rational number with denominator $m$, but other rational numbers occur for the examples in [3].

[2] The proof in Section 3 that $L$ is $G^s$ hypoelliptic for all $s \geq q/p$ applies equally well to such generalizations. The very simple analysis given in Section 4 for $s < q/p$ no longer applies when the coefficients of $L$ depend on $t$, however. The techniques of [4],[1],[2] can be adapted to this situation, but additional labor would be necessary.

Examples exhibiting a rather different behavior in the Gevrey scale have been analyzed by Yu [14]: In $\mathbb{R}^5$ with coordinates $(x_1, y_1, x_2, y_2, t)$ set $X_j = \partial_{x_j}$, $Y_1 = \partial_{y_1} + x_1^{m-1}\partial_t$, and $Y_2 = \partial_{y_2} + x_2\partial_t$ for $m \geq 2$. Then the optimal exponent for Gevrey regularity of $L = X_1^2 + X_2^2 + Y_1^2 + Y_2^2$ in any neighborhood of the origin is $s = 1$ for $m = 2$, and is $s = 2$ for all even $m > 2$.[3] Thus in this family, the optimal Gevrey exponent is independent of the degree of degeneracy of the operator.

Throughout the paper, $C, c, \varepsilon, \varepsilon', \delta$ denote strictly positive constants whose values are permitted to vary freely from one line to the next. Generally $C$ is a large constant, while the others are small.

## 2. Characterizations of Gevrey Classes

For $y = (y_1, \ldots y_n) \in \mathbb{C}^n$ define $\langle y \rangle = (1 + \Sigma y_j^2)^{1/2}$, which is well defined and holomorphic in a conic neighborhood $\Gamma$ of $\mathbb{R}^n$. Let $(x, \xi)$ be coordinates on $\mathbb{C}^n \times \mathbb{C}^n$. For each $\gamma \in [0, 1]$ define a differential form $\omega$ on a conic neighborhood of $\mathbb{R}^n$ in $\mathbb{C}^n$ by

$$\omega = dx_1 \wedge \cdots \wedge dx_n \wedge d(\xi_1 + ix_1\langle\xi\rangle^\gamma) \wedge \cdots \wedge d(\xi_n + ix_n\langle\xi\rangle^\gamma)$$

and define a function $\alpha_\gamma$ by

$$\omega = \alpha_\gamma(x, \xi) \, dx_1 \wedge \cdots \wedge dx_n \wedge d\xi_1 \wedge \cdots \wedge d\xi_n.$$

For $u \in \mathcal{E}'(\mathbb{R}^n)$, for $(x, \xi) \in \mathbb{C}^n \times \Gamma$, and for $0 \leq \gamma \leq 1$ define

$$\mathcal{F}_\gamma u(x, \xi) = \left\langle u, e^{i(x-x')\cdot\xi - \langle\xi\rangle^\gamma(x-x')^2} \alpha_\gamma(x - x', \xi) \right\rangle,$$

where the pairing is that of distributions with test functions, with respect to the variable $x'$. The exponent in this expression may be written more explicitly as

$$i\sum_{j=l}^{n}(x_j - x'_j)\xi_j - \langle\xi\rangle^\gamma \sum_{j=1}^{n}(x_j - x'_j)^2.$$

The coefficient $\alpha_\gamma$ is holomorphic with respect to $x$, and equals $1 + O(\langle\xi\rangle^{\gamma-1})$ for $x$ in any bounded subset of $\mathbb{C}^n$. Moreover $\alpha_\gamma^{-1}(x, \xi)$ is uniformly bounded, and equals $1 + O(\langle\xi\rangle^{\gamma-1})$ for large $|\xi|$, provided that $x$ is restricted to a sufficiently small but fixed neighborhood of the origin.

---

[3]Gevrey hypoellipticity for the indicated ranges of exponents follows from theorems of Treves and of Tartakoff for $m = 2$ and of Derridj and Zuily for $m > 2$; that hypoellipticity does not hold for smaller exponents was proved by Yu.

**Lemma 2.1.** *Let $0 \leq \gamma \leq 1$. For any $u \in C^{n+1}(\mathbb{R}^n)$ having compact support, for any $x \in \mathbb{R}^n$,*

$$u(x) = (2\pi)^{-n} \int_{\mathbb{R}^n} \mathcal{F}_\gamma u(x, \xi) \, d\xi.$$

*Proof.* Write the Fourier inversion formula as

$$u(x) = (2\pi)^{-n} \lim_{\varepsilon \to 0} \iint_{\mathbb{R}^n \times \mathbb{R}^n} u(x') e^{i(x-x') \cdot \xi} e^{-\varepsilon \xi^2} \, dx' \wedge d\xi$$

where $\xi^2 = \sum_j \xi_j^2$ for $\xi \in \mathbb{C}^n$. Consider this integral with some small value of $\varepsilon$, and shift the contour of integration from $\mathbb{R}^n \times \mathbb{R}^n \subset \mathbb{R}^n \times \mathbb{C}^n$ to the contour

$$\Gamma(y, \eta) = (y, \eta + i \langle \eta \rangle^\gamma (x - y)) \text{ for } (y, \eta) \in \mathbb{R}^n \times \mathbb{R}^n.$$

The integrand is holomorphic with respect to $\xi$. It also decays rapidly, for $(y, \eta) \in \mathbb{R}^{2n}$, as

$$|(x', \xi)| = |\Gamma(y, \eta, t)| \to \infty,$$

while $x' = y$ remains bounded, where $\Gamma(y, \eta, t) = (1-t)(y, \eta) + t\Gamma(y, \eta)$ for $0 \leq t \leq 1$, provided that either $\gamma < 1$, or $|x - y|$ is sufficiently small for all $y$ in the support of $u$. Indeed,

$$\mathrm{Re}(-\varepsilon \xi^2) = -\varepsilon \eta^2 + \varepsilon t^2 \langle \eta \rangle^{2\gamma} (x - y)^2$$

is bounded above by $-\varepsilon \eta^2 / 2$ for all sufficiently large $\eta$, in either of these two cases. Since $dx' \wedge d\xi = \alpha_\gamma(y, \eta) \, dy \wedge d\eta$,

$$u(x) = (2\pi)^{-n} \lim_{\varepsilon \to 0} \iint_{\mathbb{R}^n \times \mathbb{R}^n} u(y) e^{i(x-y) \cdot \eta - \langle \eta \rangle^\gamma (x-y)^2} \alpha_\gamma(x-y, \eta) e^{-\varepsilon(\eta + i\langle\eta\rangle^\gamma(x-y))^2} \, dy \, d\eta.$$

In either of the two favorable cases, fixing any sufficiently largte $\eta$ and integrating by parts $n + 1$ times with respect to $y$ yields

$$\int u(y) e^{i(x-y) \cdot \eta - \langle \eta \rangle^\gamma (x-y)^2} \alpha_\gamma(x-y, \eta) e^{-\varepsilon(\eta + i\langle\eta\rangle^\gamma(x-y))^2} \, dy = O(|\eta|^{-n-1}),$$

uniformly in $\varepsilon$ and uniformly for $x$ in any bounded set. This integral converges, as $\varepsilon \to 0$, to $\mathcal{F}_\gamma u(x, \eta)$, which is therefore $O(|\eta|^{-n-1})$. Thus $(2\pi)^{-n} \int \mathcal{F}_\gamma u(x, \xi) \, d\xi$ converges absolutely, and uniformly on compact sets, to $u(x)$.

In the case $\gamma = 1$ an additional argument is needed. By introducing a partition of unity we may assume $u$ to be supported in a ball $B(x_0, \varepsilon)$ of center $x_0$ and radius $\varepsilon$, with $\varepsilon$ chosen to be sufficiently small that the inversion formula is valid for all $x \in B(x_0, 2\varepsilon)$. Setting $v(x) = (2\pi)^{-n} \int \mathcal{F}_\gamma u(x, \xi) \, d\xi$, we thus have $v \equiv 0$ in $B(x_0, 2\varepsilon) \setminus B(x_0, \varepsilon)$. But when $\gamma = 1$, a simple calculation shows that $v$ is real analytic outside the support of $u$, so $v \equiv 0$ on the complement of $B(x_0, \varepsilon)$. Thus $v \equiv u$. □

**Lemma 2.2.** *Let $\gamma \in [0,1]$ and $u \in \mathcal{E}'(\mathbb{R}^n)$. Define*

$$u_\varepsilon(x) = (2\pi)^{-n} \int_{|\xi| \leq \varepsilon^{-1}} \mathcal{F}_\gamma u(x, \xi) \, d\xi.$$

*Then $u_\varepsilon \to u$ in $\mathcal{E}'$ as $\varepsilon \to 0$.*

*Proof.* This follows from Lemma 2.1 by duality. □

**Definition.** Let $U \subset \mathbb{R}^n$ be an open neighborhood of a point $x_0$, and $s \in (0, \infty)$. The distribution $u$ is said to belong to the Gevrey class $G^s$ at $x_0$ if there exist a neighborhood $V \subset U$ of $x_0$ such that $u \in C^\infty(V)$, and a constant $C \in \mathbb{R}^+$, such that for every $x \in V$ and every multi-index $\alpha$,

$$|\partial^\alpha u(x)| \leq C^{|\alpha|+1} |\alpha|^{s|\alpha|}.$$

Moreover $u \in G^s(U)$ if $u \in G^s$ at every point of $U$.

It is clear, by the Cauchy-Kowalevska theorem, that $G^s$ hypoellipticity cannot hold for $s < 1$, so we restrict attention henceforth to the case $s \geq 1$. The principal result of this section relates three characterizations of the classes $G^s$; two of these characterizations are formulated each in two distinct ways.

**Theorem 2.3.** *Let $x_0 \in \mathbb{R}^n, u \in \mathcal{D}'(\mathbb{R}^n)$, and $s \in [1, \infty)$. Then the following five conditions are mutually equivalent:*

(1) $u \in G^s$ at $x_0$.
(2) *For every $\gamma \geq s^{-1}$ and every distribution $v \in \mathcal{E}'(\mathbb{R}^n)$ satisfying $v \equiv u$ in some neighborhood of $x_0$, there exist $C, \delta \in \mathbb{R}^+$ and a neighborhood $V$ of $x_0$ such that*

$$|\mathcal{F}_\gamma v(x, \xi)| \leq C e^{-\delta \langle \xi \rangle^{1/s}} \quad \text{for all } (x, \xi) \in V \times \mathbb{R}^n.$$

(3) *There exist $v \in \mathcal{E}'$ agreeing with $u$ in some neighborhood of $x_0$, $\gamma \in [s^{-1}, 1]$, and a neighborhood $V$ of $x_0$ such that*

$$|\mathcal{F}_\gamma v(x, \xi)| \leq C e^{-\delta \langle \xi \rangle^{1/s}} \quad \text{for all } (x, \xi) \in V \times \mathbb{R}^n.$$

(4) *There exist an open neighborhood $U \subset \mathbb{C}^n$ of $x_0$ and $C, \delta \in \mathbb{R}^+$ such that, for each $\lambda \geq 1$, there exists a decomposition*

$$u = g_\lambda + h_\lambda \quad \text{in } U \cap \mathbb{R}^n$$

*such that $g_\lambda$ is holomorphic in $U$,*

$$|g_\lambda(z)| \leq C e^{C\lambda |\operatorname{Im}(z)|} \quad \text{for all } z \in U.$$

*and*

$$|h_\lambda(x)| \leq C e^{-\delta \lambda^{1/s}} \quad \text{for all } x \in U \cap \mathbb{R}^n.$$

(5) *There exist an open neighborhood $U \subset \mathbb{C}^n$ of $x_0$ and $C, \delta \in \mathbb{R}^+$ such that, for each $\lambda \geq 1$, there exists a decomposition*
$$u = g_\lambda + h_\lambda \quad \text{in } U \cap \mathbb{R}^n$$
*such that $g_\lambda$ is holomorphic in $\{z \in U : |\operatorname{Im}(z)| \leq \lambda^{-(s-1)/s}\} = U_\lambda$,*
$$|g_\lambda(z)| \leq C \quad \text{for all } z \in U_\lambda,$$
*and*
$$|h_\lambda(x)| \leq Ce^{-\delta \lambda^{1/s}} \quad \text{for all } x \in U \cap \mathbb{R}^n.$$

The implication (2) $\Rightarrow$ (3) is tautologous, and (4) $\Rightarrow$ (5) will be shown to be rather superficial as well.

*Proof that* (1) $\Rightarrow$ (2). Assume that $u \in G^s$ at $x_0$, $\gamma \geq s^{-1}$, and $v \in \mathcal{E}'(\mathbb{R}^n)$ agrees with $u$ in some open neighborhood $U$ of $x_0$. Fix a cube $Q$ centered at $x_0$, whose closure is contained in $U$, and fix any relatively compact open subset $V$ of the interior of $Q$. Then uniformly for $(x, \xi) \in V \times \mathbb{R}^n$,
$$\mathcal{F}_\gamma v(x, \xi) = \int_Q e^{i(x-y)\cdot\xi} \cdot \left(v(y)e^{-\langle\xi\rangle^\gamma (x-y)^2} \alpha_\gamma(x-y, \xi)\right) dy + O(e^{-\delta\langle\xi\rangle^\gamma}). \tag{2.1}$$

Assume without loss of generality that $|\xi_1| \geq n^{-1}|\xi| \geq 1$, and let $N$ be a large parameter to be chosen below. Integrating by parts $N$ times with respect to $y_1$ in (2.1) yields
$$\pm i^N |\xi_1|^{-N} \int_Q e^{i(x-y)\cdot\xi} \partial_{y_1}^N \left(v(y)e^{-\langle\xi\rangle^\gamma(x-y)^2}\alpha_\gamma(x-y, \xi)\right) dy,$$
modulo boundary terms which will be discussed below. To simplify notation set $\lambda = \langle\xi\rangle$. For $(x, y)$ in any compact set,
$$|\partial_{y_1}^k e^{-\langle\xi\rangle^\gamma (x-y)^2}| \leq C^{k+1} \lambda^{\gamma k/2} k^{k/2}$$
for any $k \geq 0$. Since $v \equiv u \in G^s$ on a neighborhood of the closure of $Q$, Leibniz's rule yields
$$|\xi_1|^{-N} \cdot |\partial_{y_1}^N \left(v(y)e^{-\langle\xi\rangle^\gamma(x-y)^2}\alpha_\gamma(x-y,\xi)\right)| \leq \lambda^{-N} C^{N+1} \max_{a+b=N} N^{sa}\lambda^{\gamma b/2} N^{b/2}.$$
where $\lambda = \langle\xi\rangle$. Fix a small constant $\varepsilon > 0$, and choose $N$ so that $|N - \varepsilon\lambda^{1/s}| \leq 1$. Then the right-hand side of the preceding display is bounded by
$$C^{\varepsilon\lambda^{1/s}+1} \max_{a+b=N} (\varepsilon\lambda^{1/s})^{sa} \lambda^{\gamma b/2} (\varepsilon\lambda^{1/s})^{b/2} \lambda^{-N}$$
$$\leq C^{\varepsilon\lambda^{1/s}+1} \max_{0 \leq b \leq N} \varepsilon^{sa+b/2} \lambda^{(\frac{b\gamma}{2}+\frac{b}{2s}-b)}.$$

Now $\gamma/2 + 1/(2s) - 1 \leq \frac{1}{2} + \frac{1}{2} - 1 = 0$, and $sa + \frac{1}{2}b \geq \frac{N}{2}$, so this last quantity is bounded above by

$$C^{\varepsilon\lambda^{1/s}+1}\varepsilon^{N/2} \leq C \cdot (C\varepsilon^{1/2})^{\varepsilon\lambda^{1/s}} \leq Ce^{-\delta\langle\xi\rangle^{1/s}}$$

for some $C, \delta \in \mathbb{R}^+$ if $\varepsilon$ is chosen to be sufficiently small.

It remains to handle the boundary terms arising from integration by parts in (2.1), which necessitates an improved bound for the derivatives of $\exp(-\langle\xi\rangle^\gamma(x-y)^2)$ away from the diagonal. Setting $f(x) = \exp(-\langle\xi\rangle^\gamma x^2)$, one has, for $x$ in any compact set disjoint from the origin,

$$|\partial_{x_1}^k f(x)| \leq C^{k+1}k^k e^{-\delta\langle\xi\rangle^\gamma},$$

as is seen by applying Cauchy's integral formula to the holomorphic function $f$ on disks disjoint from the origin. Thus for $x \in V$ and $y \in \partial Q$, for any $M \leq N$,

$$\begin{aligned}|\xi_1|^{-M} \cdot |\partial_{y_1}^M\bigl(v(y)e^{-\langle\xi\rangle^\gamma(x-y)^2}\alpha_\gamma(x-y,\xi)\bigr)| &\leq \lambda^{-M}C^{M+1}e^{-\delta\langle\xi\rangle^\gamma}\max_{a+b=M} M^{sa}M^b \\ &\leq C^{M+1}M^{sM}\lambda^{-M}e^{-\delta\lambda^{1/s}} \\ &\leq C^{M+1}(N^s/\lambda)^M e^{-\delta\lambda^{1/s}} \\ &\leq C(C\varepsilon^s)^M e^{-\delta\lambda^{1/s}} \\ &\leq Ce^{-\delta\lambda^{1/s}}.\end{aligned}$$

Thus each boundary term is $O(\exp(-\delta\lambda^{1/s}))$, and there are $2N \leq C\lambda^{1/s}$ of them. $\square$

The transform $\mathcal{F}_\gamma v(z,\xi)$ extends, for each $\xi \in \mathbb{R}^n$, to an entire holomorphic function of $z \in \mathbb{C}^n$. Observe that the same reasoning as above gives

$$|\mathcal{F}_\gamma v(z,\xi)| \leq Ce^{-\delta\langle\xi\rangle^{1/s}}e^{C\langle\xi\rangle|\text{Im}(z)|} \qquad (2.2)$$

for $z$ in a sufficiently small neighborhood $V \subset \mathbb{C}^n$ of $x_0$, provided that $v \in \mathcal{E}'$ belongs to $G^s$ at $x_0$ and that $1 \geq \gamma \geq s^{-1}$.

*Proof that* $(1) \Rightarrow (4)$. Suppose that $u \in G^s$ at $x_0$. Fix $v \in \mathcal{E}'$ that agrees with $u$ near $x_0$, and set $\gamma = s^{-1}$, so that (2.2) holds. Thus it suffices to show that (2.2) implies (4). For each large $\lambda \in \mathbb{R}^+$ define

$$g_\lambda(z) = (2\pi)^{-n}\int_{|\xi|\leq\lambda}\mathcal{F}_\gamma v(z,\xi)\,d\xi.$$

Since $|\mathcal{F}_\gamma v(z,\xi)| \leq C\exp(-\delta\langle\xi\rangle^\gamma + C\langle\xi\rangle|\text{Im}(z)|)$ in a complex neighborhood $U \subset \mathbb{C}^n$ independent of $\xi$, and since $z \to \mathcal{F}_\gamma v(z,\xi)$ is holomorphic in $U$, each $g_\lambda$ is holomorphic

on $U$ and satisfies
$$|g_\lambda(z)| \le C \int_{r=0}^{\lambda} e^{-\delta r^\gamma} e^{Cr|\text{Im}(z)|} r^{n-1} dr \le Ce^{C\lambda|\text{Im}(z)|}.$$

Moreover, for $x \in U \cap \mathbb{R}^n$,
$$|(u - g_\lambda)(x)| = |v - g_\lambda(x)| = (2\pi)^{-n} \left| \int_{|\xi| \ge \lambda} \mathcal{F}_\gamma v(x, \xi)\, d\xi \right|$$
$$\le C \int_{r=\lambda}^{\infty} e^{-\delta r^\gamma} r^{n-1} dr \le Ce^{-\varepsilon \lambda^\gamma}$$

for some $\varepsilon > 0$. □

*Proof that* (4) $\Rightarrow$ (5). For each $\lambda \in \mathbb{R}^+$ fix a decomposition $u = g_\lambda + h_\lambda$ having the stated properties. Then there exist a neighborhood $V \subset \mathbb{C}^n$ of $x_0$, relatively compact in $U$, and a constant $\varepsilon > 0$ such that for all $\lambda \ge 1$ and all $z \in V$ satisfying $|\text{Im}(z)| \le \varepsilon \lambda^{-(s-1)/s}$,
$$|g_\lambda(z) - g_{\lambda/2}(z)| \le e^{-\varepsilon \lambda^{1/s}}. \tag{2.3}$$

Indeed, the function $\log|g_\lambda - g_{\lambda/2}|$ is subharmonic, is $O(C\lambda|\text{Im}(z)|)$ for $z \in U$, and is $\le -\delta \lambda^{1/s}$ for $z \in U \cap \mathbb{R}^n$, because there it equals $\log|h_\lambda - h_{\lambda/2}|$. A simple comparison argument, applied in complex one dimensional rectangles each having one boundary segment in $\mathbb{R}^n$, then yields (2.3).

Given any large $\lambda \in \mathbb{R}^+$, choose $k$ so that $2^k \le \lambda < 2^{k+1}$ and write
$$g_\lambda = g_{\lambda/2^k} + \sum_{j=1}^{k} \left( g_{\lambda/2^{j-1}} - g_{\lambda/2^j} \right).$$

By (2.3), $g_\lambda(z)$ is $O(1)$ for all $z \in V$ satisfying $|\text{Im}(z)| \le \varepsilon' \lambda^{-(s-1)/s}$. Thus $g_\lambda$ has all the properties required in (5) of $g_\mu$, where $\mu^{-(s-1)/s} = \varepsilon' \lambda^{-(s-1)/s}$. □

*Proof that* (5) $\Rightarrow$ (2). Suppose that $v \equiv u$ in $U_0$, that $v$ is supported in $U_1$ and that $U_0$ is relatively compact in $U_1$. Assume that $\gamma \ge s^{-1}$. To estimate $\mathcal{F}_\gamma v(x, \xi)$ for $x \in U_0$, assume without loss of generality that $|\xi_1| \ge |\xi|/n \ge 1$ and that $\xi_1 > 0$. Fix $\eta \in C_0^\infty(U_1)$ satisfying $\eta \equiv 1$ on a neighborhood of the closure of $U_0$. Then for $x \in U_0$,
$$\mathcal{F}_\gamma v(x, \xi) = \mathcal{F}_\gamma (g_\lambda \eta)(x, \xi) + O(e^{-\delta \lambda^{1/s}})$$

where $\lambda$ is defined to be $|\xi|$, and $u = g_\lambda + h_\lambda$ in $U_1$ with $g_\lambda, h_\lambda$ as in (5). Translate the coordinate system to have origin at $x_0$, and choose $U_0$ of the form
$$U_0 = \{x : |x_j| < r \text{ for all } j\},$$

where $r$ is small enough that $x \in U_1$ whenever each $|x_j| < 3r$. Choose $\eta \in C_0^\infty(U_1)$ so $\eta(x) = 1$ whenever each $|x_j| \leq 2r$. Fix $\varphi \in C^1(\mathbb{R})$ such that $0 \leq \varphi \leq \frac{1}{2}$, $\varphi(t) > 0$ for $|t| \leq 2r$, and $\varphi(t) = 0$ for all $|t| \geq 3r$. Letting $\varepsilon > 0$ be small, shift the contour of integration in the definition of $\mathcal{F}_\gamma(g_\lambda \eta)(x, \xi)$ to

$$\phi(y) = \left(y_1 + i\varepsilon\lambda^{-(s-1)/s}\varphi(y_1), y_2, \ldots, y_n\right)$$

to obtain

$$\mathcal{F}_\gamma(g_\lambda \eta)(x, \xi)$$
$$= c\int_{\mathbb{R}^n} e^{i(x-\phi(y))\cdot\xi - \langle\xi\rangle^\gamma(x-\phi(y))^2} \alpha_\gamma(x-\phi(y), \xi)\left(1 + i\varepsilon\lambda^{-(s-1)/s}\varphi'(y_1)\right) \cdot g_\lambda\eta(\phi(y))\, dy.$$

Since $\eta(y) \equiv 1$ for all $y \in \mathbb{R}^n$ such that $\phi(y) \neq y$, we interpret $\eta(\phi(y))$ to be $\eta(y)$. The function $g_\lambda$ is holomorphic in a region containing $\phi(y)$ for all $y$ in the support of $\eta$, so $g_\lambda(\phi(y))$ is defined and the change of contour is justified.

Now

$$\operatorname{Re}\left(i(x-\phi(y))\cdot\xi - \langle\xi\rangle^\gamma(x-\phi(y))^2\right)$$
$$= -\langle\xi\rangle^\gamma(x-y)^2 - \xi_1\lambda^{-(s-1)/s}\varepsilon\varphi(y_1) + \langle\xi\rangle^\gamma\lambda^{-2(s-1)/s}\varepsilon^2\varphi^2(y_1)$$
$$\leq -\delta\lambda^\gamma(x-y)^2 - \delta\lambda^{1/s}\varphi(y_1)$$

for some $\delta > 0$, for all sufficiently large $\lambda = |\xi|$, if $s > 1$ or if $\varepsilon$ is chosen to be sufficiently small. The integrand is then $O(\exp(-\delta\lambda^{1/s}\varphi(y_1))) = O(\exp(-\delta'|\xi|^{1/s}))$ for $|y_1| \leq 2r$, since $g_\lambda(y) = O(1)$ where $|\operatorname{Im}(y)| \leq C|\xi|^{-(s-1)/s}$, and for $|y_1| \geq 2r$ is

$$O\left(\exp(-\delta\lambda^\gamma |x_1 - y_1|^2)\right) = O\left(\exp(-c|\xi|^{1/s})\right)$$

if $x \in U_0$. $\square$

*Proof that* (3) $\Rightarrow$ (1). Suppose that $v \equiv u$ near $x_0$, $v \in \mathcal{E}'$, $\gamma \geq s^{-1}$, and

$$\mathcal{F}_\gamma v(x, \xi) = O(\exp(-c\langle\xi\rangle^{1/s})) \quad \text{for all } x \in \mathbb{R}^n$$

in a fixed neighborhood of $x_0$. There is also the trivial bound

$$\mathcal{F}_\gamma v(z, \xi) = O(\exp(C\langle\xi\rangle)),$$

valid for all $z$ in a fixed complex neighborhood of $x_0$. Subharmonicity of $\log|\mathcal{F}_\gamma v(\cdot, \xi)|$ and a simple comparison argument based on these two bounds yields

$$|\mathcal{F}_\gamma v(z, \xi)| \leq Ce^{-\delta\langle\xi\rangle^{1/s} + C|\operatorname{Im}(z)|\langle\xi\rangle},$$

valid for all $\xi \in \mathbb{R}^n$ and all $z \in \mathbb{C}^n$ in a smaller neighborhood of $x_0$. In particular, there exist $C, \delta \in \mathbb{R}^+$ such that

$$|\mathcal{F}_\gamma v(z, \xi)| \leq C e^{-\delta \langle \xi \rangle^{1/s}}$$

whenever $|z| < \delta$ and $|\text{Im}(z)| < \delta \langle \xi \rangle^{-(s-1)/s}$. Therefore for all $x \in \mathbb{R}^n$ satisfying $|x| < \delta/2$, we may apply Cauchy's inequalities to $\mathcal{F}_\gamma v(\cdot, \xi)$ in the polydisk of radii $\delta' \langle \xi \rangle^{-(s-1)/s}$ centered at $x$ to obtain

$$|\partial_x^\alpha \mathcal{F}_\gamma v(x, \xi)| \leq C^{|\alpha|+1} |\alpha|^{|\alpha|} \langle \xi \rangle^{|\alpha|(s-1)/s} e^{-\delta \langle \xi \rangle^{1/s}}$$

for all multi-indices $\alpha$ and all $\xi \in \mathbb{R}^n$. The inversion formula for $\mathcal{F}_\gamma$ and these bounds together imply that $u = v \in C^\infty$ in the region where $|x| < \delta/2$, and

$$
\begin{aligned}
|\partial_x^\alpha u(x)| &\leq C^{|\alpha|+1} |\alpha|^{|\alpha|} \int_{\mathbb{R}^n} \langle \xi \rangle^{|\alpha|(s-1)/s} e^{-\delta \langle \xi \rangle^{1/s}} \, d\xi \\
&\leq C^{|\alpha|+1} |\alpha|^{|\alpha|} \int_0^\infty t^{|\alpha|(s-1)/s} e^{-\delta t^{1/s}} t^{n-1} \, dt \\
&\leq C^{|\alpha|+1} |\alpha|^{|\alpha|} \int_0^\infty r^{|\alpha|(s-1)} e^{-\delta r} r^{ns-1} \, dr \\
&= C^{|\alpha|+1} |\alpha|^{|\alpha|} \Gamma(|\alpha|(s-1) + ns) \\
&\leq C^{|\alpha|+1} |\alpha|^{s|\alpha|}.
\end{aligned}
$$

Thus $u \in G^s$ at $x_0$. $\square$

We have proved the chain of implications $(1) \Rightarrow (2.2) \Rightarrow (4) \Rightarrow (5) \Rightarrow (2) \Rightarrow (3) \Rightarrow (1)$, so the proof of the theorem is complete.

## 3. Gevrey Regularity

Let $1 \leq p \leq q$ be integers and set

$$L = \partial_x^2 + (x^{p-1} \partial_{t_1})^2 + (x^{q-1} \partial_{t_2})^2, \tag{3.1}$$

with coordinates $(x, t) \in \mathbb{R} \times \mathbb{R}^2$. Let $\eta = (\xi, \tau)$ be dual coordinates. Our aim is to prove $L$ to be $G^s$ hypoelliptic for all $s \geq q/p$. The proof is based on characterizations (2) and (5) of $G^s$, together with the following lemma. To simplify notation, set

$$E = E(t') = e^{i(t-t') \cdot \tau - \langle \eta \rangle^{p/q} (t-t')^2}.$$

**Lemma 3.1.** *Let $L$ take the form (3.1), for some $1 \leq p \leq q$. There exist $c, c', \delta \in \mathbb{R}^+$ and an open neighborhood $U \subset \mathbb{R} \times \mathbb{C}^2$ of the origin such that for each $\eta = (\xi, \tau) \in \mathbb{R}^3$* *satisfying $|\xi| \leq |\tau|$ and each $(x,t) \in U$ there exists $g \in C^0(U)$, holomorphic with respect to $t$, satisfying*

$$L(Eg)(x', t') = \alpha_{p/q}(x - x', t - t', \eta) e^{i(x-x')\cdot\xi - \langle\eta\rangle^{p/q}(x-x')^2} E(x', t') + O(e^{-\delta\langle\eta\rangle^{p/q}})$$

*in $U \cap \mathbb{R}^3$,*

$$g = O(1) \text{ in } U,$$

*and for $(x', t') \in U \cap \mathbb{R}^3$,*

$$g(x', t') = O(e^{-\delta|\tau|^{p/q}}) \quad \text{if } |x'| > c \text{ and } |x| < c',$$

*uniformly in $(x, t, \eta)$.*

The method of proof is the same as that of Lemmas 5.1 of [1] and 4.1 of [2]. Before presenting the details, we indicate how it implies Gevrey hypoellipticity. Set $\gamma = p/q$.

Suppose that $W \subset \mathbb{R}^3$ is open, $u \in \mathcal{D}'(W)$, and $Lu \in G^s(W)$ for some $s \geq \gamma^{-1} = q/p$. Because $L$ is elliptic wherever $x \neq 0$, $u \in G^s$ on $W \setminus \{(x,t) : x = 0\}$ [7]. Since $L$ is invariant under translation with respect to $t$, it suffices to show that if $0 \in W$, then $u \in G^s$ at 0. Fix a relatively compact open set $V \subset U \cap W$ and $h \in C_0^\infty(U \cap W)$ satisfying $h \equiv 1$ in a neighborhood of $\overline{V}$. Replace $u$ by $u \cdot h$. Then $u, Lu \in C^\infty(U \cap W)$, because $L$ is $C^\infty$ hypoelliptic.

Consider any $(x,t) \in V$ and $\eta = (\xi, \tau) \in \mathbb{R} \times \mathbb{R}^2$ with $|\xi| \leq |\tau|$, and let $g = g_{(x,t,\eta)}$ satisfy the conclusions of Lemma 3.1. Then

$$\mathcal{F}_\gamma u(x, t, \eta) = \langle u, L(Eg)\rangle + O(e^{-\delta\langle\eta\rangle^\gamma})$$
$$= \langle Lu, Eg\rangle + O(e^{-\delta\langle\eta\rangle^\gamma}).$$

By the characterization (5) of $G^s$, there exists an open set $0 \in V_1 \subset \mathbb{C}^3$ such that for each $\eta$, $Lu$ may be decomposed in $V_1 \cap \mathbb{R}^3$ as

$$Lu = G + O(e^{-\delta\langle\eta\rangle^{1/s}}),$$

where $G$ is holomorphic and $O(1)$ in

$$\{(x,t) \in V_1 : |\text{Im}(x,t)| < \langle\eta\rangle^{-(s-1)/s}\},$$

uniformly in $\eta$. For all $(x,t)$ in a compact subset of an open ball $B \subset V \cap V_1$ centered at 0, then,

$$\mathcal{F}_\gamma u(x, t, \eta) = \int_B Lu(x', t') g(x', t') e^{i(t-t')\cdot\tau - \langle\eta\rangle^\gamma(t-t')^2} dx' dt' + O(e^{-\delta\langle\eta\rangle^\gamma})$$
$$= \int_B G(x', t') g(x', t') e^{i(t-t')\cdot\tau - \langle\eta\rangle^\gamma(t-t')^2} dx' dt' + O(e^{-\delta\langle\eta\rangle^{1/s}}).$$

The hypothesis $s \geq q/p = 1/\gamma$ has been invoked in passing from the first line to the second, while the fact that $g(x', t') \exp(-\langle \eta \rangle^\gamma (t - t')^2)$ is $O(\exp(-\delta|\tau|))$ for $(x', t')$ outside $B$ has been used to justify restricting the integral to $B$ in the first line; this bound for $g$ holds, by Lemma 3.1, provided that $B$ is chosen to be sufficiently small and $x$ is restricted to a compact subset of $B$. This last integral is $O(\exp(-\delta\langle \eta \rangle^{1/s}))$, as is seen by deforming the contour of integration with respect to $t'$ as in the proof that $(5) \Rightarrow (2)$ in Section 2, exploiting the holomorphic extendibility of both $G$ and $g$. Thus we have shown that $\mathcal{F}_\gamma u(x, t, \eta)$ is $O(\exp(-\delta\langle \eta \rangle^{1/s}))$, in the region where $|\xi| \leq |\tau|$.

There remains the elliptic region, where $|\xi| \geq |\tau|$. Then $\mathcal{F}_\gamma u(x, t, \eta)$ is analyzed in the same way, using the following simpler variant of Lemma 3.1. Let

$$E = e^{i(x-x', t-t') \cdot \eta - \langle \eta \rangle^{p/q}(x-x', t-t')^2}.$$

**Lemma 3.2.** *There exist $\delta \in \mathbb{R}^+$ and an open set $0 \in U \subset \mathbb{C}^3$ such that for each $(x, t) \in U \cap \mathbb{R}^3$ and each $\eta = (\xi, \tau) \in \mathbb{R}^3$ satisfying $|\xi| \geq |\tau|$, there exists $g$ holomorphic in $U$ satisfying*

$$L(Eg)(x', t') = \alpha_{p/q}(x - x', t - t', \eta) \cdot E(x', t') + O(e^{-\delta \langle \eta \rangle}) \quad \text{in } \mathbb{R}^3 \cap U$$

*and $g = O(1)$ in $U$, uniformly in $(x, t, \eta)$.*

This completes the proof of Gevrey regularity for all $s \geq q/p$, modulo the proofs of Lemmas 3.1 and 3.2.

To begin the proof of Lemma 3.1, replace $(x, t)$ and $(x', t')$ in the statement of the lemma by $(\tilde{x}, \tilde{t})$ and $(x, t)$, respectively, so that $L$ acts with respect to $(x, t)$, and the parameters are $(\tilde{x}, \tilde{t})$ and $\eta = (\xi, \tau)$. Define $\gamma = p/q$ and set

$$E(t) = e^{i(\tilde{t}-t) \cdot \tau - \langle \eta \rangle^\gamma (\tilde{t}-t)^2} \quad \text{and} \quad L_\eta = E^{-1} \circ L \circ E.$$

Both $E$ and $L_\eta$ depend on the parameters $\tilde{t}, \eta$. Put

$$\mathcal{A}_\tau = \partial_x^2 - x^{2(p-1)}\tau_1^2 - x^{2(q-1)}\tau_2^2.$$

The operator $\mathcal{A}_\tau$ will be regarded sometimes as an ordinary differential operator acting on functions of $x$ and depending on parameters $t, \tilde{t}, \eta$, and at other times as acting on functions of $(x, t)$. For $x \in \mathbb{R}$ and $0 \neq \tau \in \mathbb{R}^2$ set

$$w(x, t) = (\tau_1^{2/p} + \tau_1^2 x^{2(p-1)} + \tau_2^{2/q} + \tau_2^2 x^{2(q-1)})^{1/2}.$$

Fix $v \in C^\infty(\mathbb{R})$, real-valued and nonnegative, satisfying $v \equiv 0$ in a small neighborhood of 0, and $v(x) \equiv 1$ for all $|x| \geq 1$. Let $\rho > 0$ be a small parameter, which will eventually depend on $v$, and define weighted Sobolev spaces $\mathcal{H}^k_\tau(\mathbb{R})$ with norms

$$\|f\|^2_{\mathcal{H}^0_\tau(\mathbb{R})} = \int_{\mathbb{R}} |f(x)|^2 w(x,\tau)^{-2} e^{\rho|\tau|^{p/q} v(x)}\, dx,$$

$$\|f\|^2_{\mathcal{H}^1_\tau(\mathbb{R})} = \int_{\mathbb{R}} [w(x,\tau)^{-2}|\partial_x f|^2 + |f|^2] e^{\rho|\tau|^{p/q} v(x)}\, dx,$$

$$\|f\|^2_{\mathcal{H}^2_\tau(\mathbb{R})} = \int_{\mathbb{R}} [w(x,\tau)^{-2}|\partial_x^2 f|^2 + |\partial_x f|^2 + w(x,\tau)^2 |f|^2] e^{\rho|\tau|^{p/q} v(x)}\, dx.$$

**Lemma 3.3.** *For all sufficiently small $|\rho|$, $\mathcal{A}_\tau : C_0^2(\mathbb{R}) \to C_0^0(\mathbb{R})$ extends to an invertible operator from $\mathcal{H}^2_\tau(\mathbb{R})$ to $\mathcal{H}^0_\tau(\mathbb{R})$, uniformly in $\rho, \tau$ for all $|\tau| \geq 1$.*

The lemma remains true if the exponent $p/q$ in the definition of the norms is replaced by any exponent belonging to $[0,1]$, but our application requires $p/q$.

*Proof.* The definitions are set up so as to imply directly that $\mathcal{A}_\tau : \mathcal{H}^2_\tau(\mathbb{R}) \to \mathcal{H}^0_\tau(\mathbb{R})$ is bounded, uniformly in $\rho, \tau$. To prove invertibility consider first the case where $\rho = 0$. For any real valued $f \in C_0^2(\mathbb{R})$,

$$-\int \mathcal{A}_\tau f \cdot f\, dx = \|\partial_x f\|^2 + \int |f|^2 (\tau_1^2 x^{2(p-1)} + \tau_2^2 x^{2(q-1)})\, dx.$$

Recall that for any $f \in C_0^2$ and any $\lambda \in \mathbb{R}^+$,

$$\lambda^{2/m} \|f\|^2 \leq C \|\partial_x f\|^2 + C \int_{\mathbb{R}} |f|^2 \lambda^2 x^{2(m-1)};$$

the case $\lambda = 1$ is elementary, and the general case then follows by scaling. Combining this with the preceding inequality yields

$$\|\partial_x f\|^2 + \int |f|^2 w(x,\tau)^2 \leq -C \int \mathcal{A}_\tau f \cdot f\, dx \leq C \|\mathcal{A}_\tau f\|_{\mathcal{H}^0_\tau} \cdot \left( \int |f|^2 w(x,t)^2 \right)^{1/2}$$

so that

$$\|\partial_x f\|^2 + \int |f|^2 w(x,\tau)^2 \leq C \|\mathcal{A}_\tau^2 f\|_{\mathcal{H}^0_\tau}.$$

We have thus controlled two of the three terms in the definition of the $\mathcal{H}^2_\tau$ norm. From the identity $\partial_x^2 = \mathcal{A}_\tau + \tau_1^2 x^{2(p-1)} + \tau_2^2 x^{2(q-1)}$, the a priori inequality

$$\|f\|^2_{\mathcal{H}^2_\tau(\mathbb{R})} \leq C \|\mathcal{A}_\tau f\|^2_{\mathcal{H}^0_\tau(\mathbb{R})} \quad \text{for all } f \in C_0^2 \tag{3.2}$$

thus follows, uniformly in $\tau$, for $\rho = 0$. The corresponding inequality for $|\rho|$ small[4] and $|\tau| \geq 1$ follows by conjugation with $\exp(\rho|\tau|^\gamma v(x)/2)$, using the fact that $w(x,\tau) \geq |\tau| \geq |\tau|^\gamma$ on the support of $\nabla v$. Now $\mathcal{H}^2_\tau(\mathbb{R})$ embeds compactly in $\mathcal{H}^0_\tau(\mathbb{R})$, and $\mathcal{A}_\tau$

---

[4]When $p/q$ is strictly less than 1, $|\rho|$ need not be chosen to be small.

is symmetric in $L^2(\mathbb{R}, dx)$; the invertibility follows from these two facts and from (3.2) as in the proof of Lemma 3.1 of [1]. Uniformity follows from the uniformity of (3.2). □

For any open disk $D \subset \mathbb{C}^2$ centered at 0, define Sobolev spaces $\mathcal{H}_\tau^k(\mathbb{R} \times D)$ of measurable, locally square-integrable functions of $(x, t) \in \mathbb{R} \times D$ that are holomorphic with respect to $t \in D$ (for almost every $x$), for which the norms are

$$\|f\|_{\mathcal{H}_\tau^0(\mathbb{R}\times D)}^2 = \int_{\mathbb{R}\times D} |f(x,t)|^2 w(x,\tau)^{-2} e^{\rho|\tau|^{p/q} v(x)} \, dx \, dt d\bar{t},$$

$$\|f\|_{\mathcal{H}_\tau^1(\mathbb{R}\times D)}^2 = \int_{\mathbb{R}\times D} [w(x,\tau)^{-2}|\partial_x f|^2 + |f|^2] e^{\rho|\tau|^{p/q} v(x)} \, dx \, dt d\bar{t},$$

$$\|f\|_{\mathcal{H}_\tau^2(\mathbb{R}\times D)}^2 = \int_{\mathbb{R}\times D} [w(x,\tau)^{-2}|\partial_x^2 f|^2 + |\partial_x f|^2 + w(x,\tau)^2|f|^2] e^{\rho|\tau|^{p/q} v(x)} \, dx \, dt d\bar{t}.$$

Clearly $\mathcal{A}_\tau : \mathcal{H}_\tau^2(\mathbb{R} \times D) \to \mathcal{H}_\tau^0(\mathbb{R} \times D)$ is invertible, as a consequence of Lemma 3.3, uniformly in $\rho, \tau, D$ provided that $|\tau| \geq 1$ and $|\rho|$ is sufficiently small.

Fix $h \in C_0^\infty(\mathbb{R})$ satisfying $h \equiv 1$ in a neighborhood of 0. Define $\mathcal{E}_\eta$ by

$$\mathcal{A}_\tau + \mathcal{E}_\eta = \partial_x^2 + x^{2(p-1)}[-i\tau_1 + 2\langle\eta\rangle^\gamma(\tilde{t}_1 - t_1)]^2 + x^{2(q-1)}[-i\tau_2 + 2\langle\eta\rangle^\gamma(\tilde{t}_2 - t_2)]^2.$$

Define $\mathcal{R}_\eta$ by $L_\eta = \mathcal{A}_\tau + \mathcal{E}_\eta + \mathcal{R}_\eta$. The perturbation term $\mathcal{R}_\eta$ involves differentiation with respect to $t$, whereas $\mathcal{A}_\tau, \mathcal{E}_\eta$ do not. Let $D \supset D'$ be open disks in $\mathbb{C}^2$ centered at the origin, with distance $(D', \partial D) \geq \varepsilon$. Let $r$ be the radius of $D$.

**Lemma 3.4.** *If $|\rho|$ is sufficiently small then the perturbation terms $\mathcal{E}_\eta, \mathcal{R}_\eta$ satisfy the following bounds, uniformly for all $|\tau| \geq 1$.*

$$h(x) \cdot \mathcal{E}_\eta : \mathcal{H}_\tau^2(\mathbb{R} \times D) \to \mathcal{H}_\tau^0(\mathbb{R} \times D) \quad \text{with norm } O(r + |\tilde{t}|). \tag{3.3}$$

$$h \cdot \mathcal{R}_\eta : \mathcal{H}_\tau^2(\mathbb{R} \times D) \to \mathcal{H}_\tau^0(\mathbb{R} \times D') \quad \text{with norm } O(\varepsilon^{-2}|\tau|^{-2\gamma}) \tag{3.4}$$

*whenever distance* $(D', \partial D) \geq \varepsilon$.

*Proof.* (3.3) would follow from the inequality

$$|\tau|^\gamma(|x|^{p-1} + |x|^{q-1}) \leq Cw(x,\tau) \quad \forall \, x \in \text{support}(h), \forall \, |\tau| \geq 1. \tag{3.5}$$

To prove this note that $|x|^{q-1} \leq C|x|^{p-1}$ for all $x$ in the support of $h$. If $|\tau_1| \geq |\tau_2|$ then since $p \leq q$,

$$|\tau|^\gamma |x|^{p-1} = |\tau|^{p/q}|x|^{p-1} \leq C|\tau_1||x|^{p-1} \leq Cw(x,\tau).$$

If $|\tau_2| \geq |\tau_1|$ then

$$|\tau|^\gamma |x|^{p-1} \leq C|\tau_2|^{p/q}|x|^{p-1} = C|\tau_2|^{1/q} \cdot (|\tau_2|^{1/q}|x|)^{p-1}$$
$$\leq C|\tau_2|^{1/q}(1 + (|\tau_2|^{1/q}|x|)^{q-1}) \leq Cw(x,\tau) \ .$$

Thus (3.3) holds.

Note that by Cauchy's inequality relating the derivative of a holomorphic function at a point to the integral of the function over a disk, $\partial_{t_j}$ maps $H^2(D)$ to $H^2(D')$ boundedly, with norm $O(\varepsilon^{-1})$, where $H^2$ denotes the space of all holomorphic square-integrable functions, with respect to Lebesgue measure. Therefore

$$\partial_{t_j} : \mathcal{H}_\tau^k(\mathbb{R} \times D) \to \mathcal{H}_\tau^k(\mathbb{R} \times D') \quad \text{is } O(\varepsilon^{-1}) \text{ for each } k,$$

uniformly in $\tau$.

By (3.5), $h(x)x^{2(p-1)}$ and $h(x)x^{2(q-1)}$ map $\mathcal{H}_\tau^2(\mathbb{R})$ to $\mathcal{H}_\tau^0(\mathbb{R})$ with a bound $O(|\tau|^{-2\gamma})$. Thus

$$h(x)x^{2(p-1)}\partial_{t_1}^2 : \mathcal{H}_\tau^2(\mathbb{R} \times D) \to \mathcal{H}_\tau^0(\mathbb{R} \times D')$$

is $O(\varepsilon^{-2}|\tau|^{-2\gamma})$; the same goes for $h(x)x^{2(q-1)}\partial_{t_2}^2$. This implies (3.4). □

We turn to the proof of Lemma 3.1. Set

$$\psi(x,t) = e^{i(\tilde{x}-x)\xi - \langle\eta\rangle^\gamma(\tilde{x}-x)^2}\alpha_\gamma(\tilde{x}-x, \tilde{t}-t, \eta).$$

It suffices to solve

$$(\mathcal{A}_\tau + h\mathcal{E}_\eta + h\mathcal{R}_\eta)g(x,t) = \psi(x,t) + O(e^{-\delta\langle\eta\rangle^\gamma})$$

globally on $\mathbb{R} \times D_0$ for some polydisk $D_0 \subset \mathbb{C}^2$ centered at 0, for then the original equation is solved in the region where $h \equiv 1$.

Fix nested polydisks $D_\infty \subset D_2 \subset D_1 \subset \mathbb{C}^2$ centered at 0, such that each is relatively compact in the next. Let $\Lambda \in \mathbb{R}^+$ be a large constant to be chosen below. Assume $|\tau|$ to be large, and choose an integer $N$ such that $|N - \Lambda^{-1}|\tau|^\gamma| < 1$. Construct polydisks $D_2 \supset D_3 \supset \cdots \supset D_N = D_\infty$, centered at 0, satisfying

$$\text{distance } (D_{j+1}, \partial D_j) \geq c\Lambda|\tau|^{-\gamma},$$

where $c$ is a constant independent of $\Lambda, \tau$.

Define $g$ by the Neumann series

$$g = \sum_{j=0}^{N} (-1)^j [\mathcal{A}_\tau^{-1}(h\mathcal{E}_\eta + h\mathcal{R}_\eta)]^j \mathcal{A}_\tau^{-1}\psi.$$

If $D_1$ is chosen to be sufficiently small, but independent of $\eta$, then $\mathcal{A}_\tau^{-1} h \mathcal{E}_\eta$ maps $\mathcal{H}_\tau^2(\mathbb{R} \times D_j)$ to itself with a bound of at most $1/2$, for all $j, \tau$, while

$$\mathcal{A}_\tau^{-1} h \mathcal{R}_\eta : \mathcal{H}_\tau^2(\mathbb{R} \times D_j) \to \mathcal{H}_\tau^2(\mathbb{R} \times D_{j+1}) \quad \text{is } O(\Lambda^{-1}).$$

Therefore it is possible to choose $\Lambda$ sufficiently large that

$$[\mathcal{A}_\tau^{-1}(h\mathcal{E}_\eta + h\mathcal{R}_\eta)] : \mathcal{H}_\tau^2(\mathbb{R} \times D_j) \to \mathcal{H}_\tau^2(\mathbb{R} \times D_{j+1})$$

has norm at most $3/4$, for all $j$, for all sufficiently large $|\tau|$. Then

$$\|g\|_{\mathcal{H}_\tau^2(\mathbb{R} \times D_\infty)} \le C < \infty,$$

uniformly in $\eta, \tilde{x}, \tilde{t}$, provided that the $\mathcal{H}_\tau^0(\mathbb{R} \times D_1)$ norm of $\psi$ is bounded uniformly in $\eta, (\tilde{x}, \tilde{t})$. This is so, provided that $\tilde{x}$ belongs to a relatively compact subset of the region where $v \equiv 0$ and $\rho > 0$ is chosen to be sufficiently small (recall the weight $\exp(\rho |\tau|^\gamma v(x))$ in the definition of the $\mathcal{H}_\tau^k$ norms).

Finally,

$$(\mathcal{A}_\tau + h\mathcal{E}_\eta + h\mathcal{R}_\eta)g = \psi \pm [(h\mathcal{E}_\eta + h\mathcal{R}_\eta)\mathcal{A}_\tau^{-1}]^{N+1}\psi,$$

and

$$[(h\mathcal{E}_\eta + h\mathcal{R}_\eta)\mathcal{A}_\tau^{-1}]^{N+1}\psi = O(3/4)^{N+1} = O(e^{-\delta \langle \eta \rangle^\gamma})$$

in the $\mathcal{H}_\tau^0(\mathbb{R} \times D_\infty)$ norm, where $\delta$ depends on the choice of $\Lambda$ but is positive. This completes the proof of Lemma 3.1. $\square$

The proof of Lemma 3.2 is parallel, but is much simpler because the principal symbol of $L$ is nonzero where $|\xi| \ge |\tau|$. The conjugated operator $L_\eta$ is replaced by $E^{-1} L E$, where

$$E(x, t) = e^{i(\tilde{x}-x, \tilde{t}-t) \cdot \eta - \langle \eta \rangle^\gamma (\tilde{x}-x, \tilde{t}-t)^2}.$$

The ordinary differential operator $\mathcal{A}_\tau$ is now replaced by the operator defined by multiplication by the symbol $-\xi^2 - \tau_1^2 x^{2(p-1)} - \tau_2^2 x^{2(q-1)}$, which for $x$ in any bounded region is comparable to $\xi^2$, hence to $\langle \eta \rangle^2$, uniformly in all parameters. All derivatives with respect to $x$ or $t$ are now incorporated into $\mathcal{R}_\eta$. One works in the simpler Hilbert spaces $H^2(D_j)$, for a sequence of polydisks $D_j \subset \mathbb{C}^3$ centered at $0$, where $1 \le j \le N \approx \Lambda^{-1} \langle \eta \rangle$ and distance $(D_{j+1}, \partial D_j) \ge c\Lambda \langle \eta \rangle$. The result is a solution $g$ of $E^{-1} L E g = \psi + O(\exp[-\delta \langle \eta \rangle])$, in $H^2(D_\infty)$ norm, for a certain polydisk $D_\infty \subset \mathbb{C}^3$ containing the origin and independent of $\eta$. $\square$

## 4. Gevrey Irregularity

**Proposition 4.1.** *If $q/p > 1$ and $s < q/p$, then $L$ is not $G^s$ hypoelliptic in any neighborhood of $0$.*

The proof follows a well established method, based upon separation of variables and a reduction to certain eigenvalue problems. For the operators under consideration here this procedure is relatively elementary, first because a complete separation of variables is possible, second because the resulting eigenvalue problems are straightforward.

Consider the ordinary differential operators
$$\mathcal{L}_z = \partial_x^2 - x^{2(q-1)} + zx^{2(p-1)}, \quad z \in \mathbb{C}.$$

**Lemma 4.2.** [10] *For any $1 \le p \le q \in \mathbb{N}$ there exist $z \in \mathbb{R}$ and a Schwartz class function $f \ne 0$, defined on $\mathbb{R}$, satisfying $\mathcal{L}_z f \equiv 0$.*

**Lemma 4.3.** *If a linear partial differential operator $P$ is $G^s$ hypoelliptic in some open set $U$ containing a point $x_0$, then there exists $B < \infty$ such that for each $f \in L^2(U)$ satisfying $Pf = 0$ in $U$,*
$$|\partial^\alpha f(x_0)| \le B^{|\alpha|+1} |\alpha|^{s|\alpha|} \|f\|_{L^2(U)}$$
*for every multi-index $\alpha$.*

The proofs of Lemma 4.2 and of Lemma 4.3 are omitted. The former is in Oleĭnik [10]. The latter is a routine application of the Baire category theorem using the Banach spaces $X_B$ of all functions $f \in L^2(U)$ such that $Lf = 0$ and such that there exists $C < \infty$ such that $|\partial^\alpha f(x_0)| \le CB^{|\alpha|} |\alpha|^{s|\alpha|}$ for all $\alpha$, as in the case $s = 1$ discussed in [10].

*Proof of Proposition 4.1.* Let $z, f$ be as in Lemma 4.1. Suppose that $L$ were $G^s$ hypoelliptic in a bounded neighborhood $U$ of $0$, for some $s < q/p$. Fix a square root $w \in \mathbb{C}$ of $z$. For each large $\lambda \in \mathbb{R}^+$ set
$$F_\lambda(x, t_1, t_2,) = e^{i\lambda t_2} e^{\lambda^{p/q} w t_1} f(\lambda^{1/q} x).$$
Then $LF_\lambda \equiv 0$ in $\mathbb{R}^3$. Fix $k \in \{0, 1\}$ such that $\partial^k f/\partial x^k(0) \ne 0$. Then
$$|\partial_x^k \partial_{t_2}^N F_\lambda(0)| = \lambda^N \lambda^{k/q} |\partial_x^k f(0)|$$
for each $N \in \mathbb{N}$, for each $\lambda$. On the other hand,
$$\|F_\lambda\|_{C_0(U)} \le C e^{C\lambda^{p/q}}$$
for some $C < \infty$, since the factors $f$ and $\exp(i\lambda t_2)$ are uniformly bounded. Thus by Lemma 4.3, there exists $B < \infty$ such that
$$\lambda^N \le CB^N N^{sN} e^{C\lambda^{p/q}} \quad \text{for all } \lambda, N.$$

For each large $N$ set $\lambda = N^{q/p}$ to deduce that
$$N^{Nq/p} \le CB^N N^{sN} e^{CN} \quad \text{for all } N.$$
But such an inequality is clearly false as $N \to \infty$, if $q/p > s$. □

*E-mail address*: christ@math.ucla.edu

University of California, Los Angeles

*Current address*: Mathematical Sciences Research Institute, 1000 Centennial Drive, Berkeley, CA 94720